\documentclass{article}

\usepackage{amssymb}

\newtheorem{theorem}{Theorem}[section]

\newtheorem{proposition}[theorem]{Proposition}

\begin{document}

\title{A new family of series expansions for $1/\pi$ and a binomial identity}

\author{\ \\
J. Sesma\thanks{Email: javier@unizar.es} \\ \ \\
{\em Departamento de F\'{\i}sica Te\'{o}rica, Facultad de Ciencias,} \\
{\em 50009, Zaragoza, Spain}}

\maketitle

\begin{abstract}
A doubly infinite set of series expansion for $1/\pi$ are reported. They follow trivially from a formal expansion for the quotient of the values taken by the gamma function for two (complex) arguments differing by an integer plus one half, obtained by an alternative computation of the Wronskian of the modified Bessel functions. The same formal expansion allows to discover also a new binomial identity.
\end{abstract}

\bigskip



\bigskip

\section{Introduction}

Series expansions for $1/\pi$ are familiar from the pioneering work of Ramanujan \cite{bern}. Proofs of those expansions  and procedures to obtain additional ones have been given by Borwein and Borwein \cite{bor1,bor2,bor3}, Chudnovsky and Chudnovsky \cite{chud}, and Guillera \cite{gui1,gui2}, among others. Here we present a doubly infinite set of series expansions for $1/\pi$ that we believe are unknown. They result as particular cases of a formal expansion which we have encountered as we were dealing with an alternative procedure of computation of the Wronskian of the modified Bessel functions $I_\nu$ and $K_\nu$. On the other hand, another particular case of the same formal expansion allows to obtain an apparently new binomial identity.

Functions of a variable $z$ are considered along the paper. Since it is a merely auxiliary variable, there is no loss of generality in assuming $z$ to be positive.

To obtain the mentioned formal expansion, we use the known value of the Wronskian \cite[Eq.~9.6.15]{abra} \cite[Eq.~10.28.2]{nist}
\begin{equation}
\mathcal{W}\{K_\nu(z),I_\nu(z)\}=1/z\,,  \label{i2}
\end{equation}
the asymptotic expansion \cite[Eq.~9.7.2]{abra} \cite[Eq.~10.40.2]{nist}
\begin{equation}
K_\nu(z)\sim\frac{\pi^{1/2}}{2^{1/2}}\,e^{-z}\sum_{n=0}^\infty a_n(\nu)\,z^{-n-1/2}\,, \label{i3}
\end{equation}
and the ascending series expansion \cite[Eq.~9.6.47]{abra} \cite[Eq.~10.39.5]{nist}
\begin{equation}
I_\nu(z)=\frac{1}{2^\nu\,\Gamma(\nu+1)}\,e^{-z}\sum_{j=0}^\infty b_j(\nu)\,z^{j+\nu}\,,  \label{i4}
\end{equation}
with coefficients
\begin{equation}
a_n(\nu)=\frac{(\nu+1/2)_n\,\langle\nu-1/2\rangle_n}{n!\,2^n}  \qquad \textrm{and} \qquad b_j(\nu)=\frac{2^j\,(\nu+1/2)_j}{j!\,(2\nu+1)_j} \label{i5}
\end{equation}
in terms of the rising and falling factorials
\begin{eqnarray}
 (x)_0\equiv 1\,, & \qquad & (x)_n\equiv x(x+1)(x+2)\cdots(x+n-1)\,,  \label{i6}  \\
 \langle x\rangle_0\equiv 1\,, & \qquad & \langle x\rangle_n\equiv x(x-1)(x-2)\cdots(x-n+1)\,.  \label{i7}
\end{eqnarray}

We show in Sect.~2 a formal expansion for the quotient $\Gamma(\nu+1)/\Gamma(\nu+k+1/2)$ ($k$ integer) which stems from a peculiar computation of the Wronskian of the modified Bessel functions $I_\nu$ and $K_\nu$. A doubly infinite family of series for $1/\pi$ result from that expansion by taking $\nu=m$, a non-negative integer, as shown in Sect.~3. The same expansion, with $\nu=m+1/2$, allows us to obtain, in Sect.~4, a binomial identity.

\section{A formal expansion for the quotient of two gamma functions}

As a previous step, we recall
 a not very common representation of the exponential function, namely
\begin{equation}
  \exp(t)\sim \sum_{k=-\infty}^\infty\frac{t^{k+\delta}}{\Gamma(k+1+\delta)}\,, \qquad \delta\in \mathbb{C}\,, \quad |\arg(t)|<\pi\,,   \label{ii1}
\end{equation}
known as Heaviside's exponential series. It was introduced by Heaviside in the second volume of his {\em Electromagnetic Theory} (London, 1899), as quoted by Hardy \cite[Sects.~2.11 and 2.12]{hard}. It has been discussed by Naundorf \cite{naun}, who has applied it to find global solutions of linear differential equations of second order with two regular or irregular singular points. According to Definition 2.1 in \cite{naun}, the symbol $\sim$ in Equation (\ref{ii1}) refers to the facts that
\begin{eqnarray}
\textrm{(i)} & & \sum_{k=0}^\infty \,\frac{t^{k+\delta}}{\Gamma(k+1+\delta)}\quad \textrm{is an entire function}, \nonumber  \\
\textrm{(ii)} & & \sum_{k=-\infty}^{-1} \frac{t^{k+\delta}}{\Gamma(k\! +\! 1\! +\! \delta)}\; \textrm{is an asymptotic series for}\;  \exp(t)-\sum_{k=0}^\infty \frac{t^{k+\delta}}{\Gamma(n\! +\! 1\! +\! \delta)} \nonumber  \\
 & &\hspace{0pt}  \textrm{as} \quad t\to\infty \quad  \textrm{in the sector} \quad |\arg(t)|<\pi\,.  \nonumber
\end{eqnarray}
Notice that the asymptotic (in the just explained sense) expansion in the right-hand side of (\ref{ii1}) satisfies the differential equation $dy(t)/dt=y(t)$ and becomes the familiar convergent series expansion of the exponential function when $\delta$ takes any integer value.

\begin{proposition}
For arbitrary $\nu\in \mathbb{C}\setminus\{-1/2, -1,-3/2, -2, \ldots\}$  and  integer $k\geq 0$, the quotient  $\Gamma(\nu+1)/\Gamma(\nu+k+1/2)$ admits the formal expansion
\begin{equation}
\frac{\Gamma(\nu+1)}{\Gamma(\nu\! +\! k\! +\! 1/2)}\sim\frac{\pi^{1/2}}{2^{2\nu}}\, \sum_{n=0}^\infty\frac{(\nu\! +\! 1/2)_n\,\langle\nu\! -\! 1/2\rangle_n\,(\nu\! +\! 1/2)_{k+n}\,(k\! +\! \nu\! +\! 2n\! +\! 1/2)}{n!\,(k+n)!\,(2\nu+1)_{k+n}}\,,   \label{ii2}
\end{equation}
\end{proposition}
{\bf Proof:}
Let us write (\ref{i2}) in the form
\begin{equation}
z\,\mathcal{W}\{e^zK_\nu(z), e^zI_\nu(z)\}=\exp(2z)\,.    \label{ii3}
\end{equation}
Substitution of $K_\nu(z)$ and $I_\nu(z)$ by their respective expansions in (\ref{i3}) and (\ref{i4}) gives for the left-hand side of  (\ref{ii3}) the formal expansion
\begin{equation}
\frac{\pi^{1/2}}{2^{\nu+1/2}\,\Gamma(\nu+1)}\, \sum_{k=-\infty}^\infty c_k(\nu)\,z^{k+\nu-1/2}\,,   \label{ii4}
\end{equation}
where
\begin{equation}
c_k(\nu)\sim\sum_{n=0}^\infty a_n(\nu)\,b_{n+k}(\nu)\,(k+2n+\nu+1/2)\,,   \label{ii5}
\end{equation}
with $a_n(\nu)$ and $b_{n+k}(\nu)$ as given in (\ref{i5}).
In turn, the right-hand side of (\ref{ii3}), can be represented by the Heaviside's exponential series
\begin{equation}
\sum_{k=-\infty}^\infty \frac{2^{k+\nu-1/2}}{\Gamma(k+\nu+1/2)}\,z^{k+\nu-1/2}\,.   \label{ii6}
\end{equation}
Comparison of the expansions (\ref{ii4}) and (\ref{ii6}) of the two sides of (\ref{ii3}) allows one to write the relation
\begin{equation}
\frac{\pi^{1/2}}{2^{\nu+1/2}\,\Gamma(\nu+1)}\;c_k(\nu)= \frac{2^{k+\nu-1/2}}{\Gamma(k+\nu+1/2)}\,. \label{ii7}
\end{equation}
Simplification of this equation, after substitution of $c_k(\nu)$ by its expression as given by (\ref{ii5}) and (\ref{i5}), completes the proof.   \hspace{125pt} $\blacksquare$

{\bf Remark:}
Both infinite sequences $\{c_k(\nu)\}$ and $\{(2^{k+\nu-1/2})/\Gamma(k+\nu+1/2)\}$ obey the recurrence relation
\begin{equation}
(k+\nu-1/2)\,y_k-2\,y_{k-1}=0\,.          \label{ii8}
\end{equation}
Since this is a first order difference equation, whose solution is unique up to a multiplicative constant, those sequences must be proportional. Our proposition unveils the proportionality constant, as given in (\ref{ii7}). Nevertheless, the convergence of the expansion in the right-hand side of  (\ref{ii5}) has not been proved and (\ref{ii2}) should be seen as purely formal. Numerical exploration allows to conjecture that it is an asymptotic expansion of the left-hand side, as a function of $k$, for $k\to\infty$. It seems to become useful, from the computational point of view, for $k\gtrsim 10$.

\section{Expansions for $1/\pi$}

In the particular case of being $\nu=m$, a non-negative integer, it is not difficult to see that the resulting series in the right-hand side of (\ref{ii2}),
\begin{equation}
 \sum_{n=0}^\infty
\frac{(m\! +\! 1/2)_n\,\langle m\! -\! 1/2\rangle_n\,(m\! +\! 1/2)_{k+n}\,(k\! +\! m\! +\! 2n\! +\! 1/2)}{n!\,(k+n)!\,(2m+1)_{k+n}}\,,   \label{iii1}
\end{equation}
turns out to be convergent provided $k\geq 2$. In fact, for $n>m$, the successive terms alternate in sign, decrease monotonously in absolute value, and go to 0 as $n\to\infty$ (Leibniz's test). Therefore, one is allowed to write, for $k\geq 2$,
\begin{equation}
\frac{\Gamma(m+1)}{\Gamma(m\! +\! k\! +\! 1/2)} = \frac{\pi^{1/2}}{2^{2m}}\,  \sum_{n=0}^\infty
\frac{(m\! +\! 1/2)_n\,\langle m\! -\! 1/2\rangle_n\,(m\! +\! 1/2)_{k+n}\,(k\! +\! m\! +\! 2n\! +\! 1/2)}{n!\,(k+n)!\,(2m+1)_{k+n}}\,.   \label{iii2}
\end{equation}
Replacement of the gamma function by their values and multiplication of both sides of this equation by $(1/2)_{m+k}\pi^{-1/2}/m!$ leads to  the family of expansions for $1/\pi$
\begin{equation}
\frac{1}{\pi}=\frac{(1/2)_{m+k}}{2^{2m}}\,\frac{(2m)!}{m!}\, \sum_{n=0}^\infty
\frac{(m\! +\! 1/2)_n\,\langle m\! -\! 1/2\rangle_n\,(m\! +\! 1/2)_{k+n}\,(k\! +\! m\! +\! 2n\! +\! 1/2)}{n!\,(k+n)!\,(2m+k+n)!}\,,   \label{iii3}
\end{equation}
with $m=0, 1, 2, \ldots$ and $k=2, 3, 4, \ldots$.
In the particular case of $m=0$, the resulting sub-family is
\begin{equation}
\frac{1}{\pi}=(1/2)_{k}\,\sum_{n=0}^\infty (-1)^n\,\left[
\frac{(1/2)_n}{n!}\right]^3\,\frac{(n\! +\! 1/2)_k\,(k\! +\! 2n\! +\! 1/2)}{\left[(n+1)_k\right]^2}\,, \quad k=2,3,4,\ldots\,,  \label{iii4}
\end{equation}
expansions which resemble those reported in \cite[Eqs.~(119) to (122)]{weis}. Given the structure of (\ref{iii4}), namely
\begin{equation}
\frac{1}{\pi}=\sum_{n=0}^\infty (-1)^n\left[\frac{(1/2)_n}{n!}\right]^3\! f_k(n), \,\; \textrm{with} \,\; f_k(n)\equiv (1/2)_{k}\, \frac{(n\! +\! 1/2)_k\,(k\! +\! 2n\! +\! 1/2)}{\left[(n+1)_k\right]^2}\,,  \label{iii5}
\end{equation}
additional series expansions for $1/\pi$ of the form
\begin{equation}
\frac{1}{\pi}=\sum_{n=0}^\infty (-1)^n\,\left[\frac{(1/2)_n}{n!}\right]^3\,g(n)  \label{iii6}
\end{equation}
can be written if one takes for $g(n)$  a ``normalized" linear combination of several $f_k(n)$ arbitrarily chosen,
\begin{equation}
g(n) = \frac{\sum_k\, \alpha_k\,f_k(n)}{\sum_k\, \alpha_k}\,,  \label{iii7}
\end{equation}
with arbitrary (even complex) coefficients $\alpha_k$.

\section{A binomial identity}

All the formal expansions in Sect.~2 become convergent series or finite sums in the particular
case of being $\nu=m\! +\! 1/2$, with $m\in \{0, \mathbb{N}\}$. Since $a_n(m\! +\! 1/2)=0$ for $n>m$, the sum in the right-hand side of (\ref{ii5}) becomes finite. Consequently, the same is true for that in the right-hand side of (\ref{ii2}), which adopts the form
\begin{equation}
\frac{\Gamma(m+3/2)}{\Gamma(m+k+1)}=\frac{\pi^{1/2}}{2^{2m+1}}\, \sum_{n=0}^m
\frac{(m+1)_n\,\langle m\rangle_n\,(m+1)_{k+n}\,(k+m+2n+1)}{n!\,(k+n)!\,(2m+2)_{k+n}}\,.   \label{iv1}
\end{equation}
Writing the gamma functions and the rising and falling factorials in this equation in terms of ordinary factorials, and grouping these in binomial coefficients, one obtains the binomial identity
\begin{equation}
\sum_{n=0}^m {m \choose n}{m+n+k \choose m}{2m+n+k \choose n+m}^{-1}\frac{m+2n+k+1}{2m+n+k+1}=1   \label{iv2}
\end{equation}
for arbitrary non-negative integers $m$ and $k$. Obviously, it can be written also in the form
\begin{equation}
\sum_{n=0}^m {m \choose n}{n+k \choose m}{m+n+k \choose n+m}^{-1}\frac{2n+k+1}{n+m+k+1}=1\,, \quad k\geq m \in \{0,\mathbb{N}\}\,.  \label{iv3}
\end{equation}
This identity seems to be new. We have searched for it in several specialized publications \cite{boya,goul,prud,rior} but we have not been able to find it.

\section*{Acknowledgments}
This work was supported by Gobierno de Arag\'on (Project 226223/1) and Mi\-nis\-te\-rio de Ciencia e Innovaci\'on (Project MTM215-64166).

\end{document}